\documentclass[12pt]{article}

\usepackage{psfrag}
\usepackage{amssymb,amsmath}
\usepackage{texdraw}
\usepackage{graphics}
\usepackage{graphicx}
\usepackage{color}
\usepackage{authblk}

\def\div{\mathrm{div}}

\def\R{\mathbb{R}}

\def\N{\mathbb{N}}

\def\dsp{\displaystyle}

\def\div{{\rm div}}

\newtheorem{thm}{Theorem}
\newtheorem{lem}{Lemma}
\newtheorem{cor}[thm]{Corollary}

\newtheorem{rem}{Remark}

\def\0{{\bf 0}}
\def\D{\mathbb D}
\def\f{\mathbf f}
\def\n{\mathbf n}
\def\N{\mathbb N}
\def\R{\mathbb{R}}
\def\s{\mathbf s}
\def\u{\mathbf u}
\def\v{\mathbf v}
\def\w{\mathbf w}
\def\W{\mathbb W}

\def\bsigma{\boldsymbol{\sigma}}
\def\btau{\boldsymbol{\tau}}

\def\div{\mathrm{div}\,}

\def\Re{\mathrm{Re}}
\def\Di{\mathrm{D}}
\def\We{\mathrm{We}}

\begin{document}




\title{Mathematical analysis of the stationary Oldroyd model with diffusive stress}


\author{Laurent Chupin\thanks{Universit\'e Blaise Pascal, Laboratoire de Math\'ematiques (CNRS UMR 6620),\\ Campus des C\'ezeaux, 63177 Aubi\`ere cedex, France,\\ 
\texttt{laurent.chupin@math.univ-bpclermont.fr}}\ \ and S\'ebastien Martin\thanks{Universit\'e Paris Descartes, Laboratoire MAP5 (CNRS UMR 8145),\\ 45 rue des Saints-P\`eres, 75270 Paris cedex 06, France,\\ 
\texttt{sebastien.martin@parisdescartes.fr}}}
\maketitle

\begin{abstract}
We present the mathematical analysis of the stationary Oldroyd model with diffusive stress: existence and uniqueness of weak solutions is shown if the source terms are small enough or if the Reynolds and Weissenberg numbers are small enough. Besides, in the corotational model, this condition on the data can be relaxed for the existence result. Finally, strong solutions are obtained with additional regularity on the data.
\end{abstract}

{\bf Keywords.} viscoelastic fluid, Oldroyd model, diffusive stress\medskip

{\bf MSC.} 76D03, 76A10, 35Q35


\section{Introduction}

The Oldroyd model describes the behaviour of a viscoelastic fluid. Its principle is built upon a description of the shear stress that interpolates between a purely viscous contribution and a purely elastic contribution. Let~$\Omega$ be a bounded open set in~$\R^3$ and~$\f$ be a given vector function on~$\Omega$.
We are looking for a vector function~$\u : \Omega \mapsto \R^3$, a scalar function~$p:\Omega \mapsto \R$ and a symmetric tensor function~$\bsigma:\Omega \mapsto \R^{3 \times 3}$, representing the velocity, the pressure and the elastic extra-stress of the fluid satisfying the following set of equations:
\begin{equation}\label{problem-strong}
\left\{
\begin{aligned}
& \Re \, ( \u\cdot \nabla \u ) - (1-r)\Delta \u + \nabla p = \div \bsigma + \f, & \text{in~$\Omega$},\\
& \div \u = 0, & \text{in~$\Omega$},\\
& \We \, ( \u \cdot \nabla \bsigma + g_{a}(\nabla \u , \bsigma) ) + \bsigma - \Di\, \Delta \bsigma = 2r \D(\u), & \text{in~$\Omega$},\\
& \u = \0, & \text{on~$\partial \Omega$},\\
& \Di\, \partial_\n \bsigma = \0, & \text{on~$\partial \Omega$},
\end{aligned}
\right.
\end{equation}
\noindent Here, $\mathrm{Re}$ denotes the Reynolds number that quantifies the inertial effects in the fluid flow; $\mathrm{We}$ is the Weissenberg number related to a relaxation time that characterizes the elasticity of the fluid; $\Di$ is a diffusive parameter of the elastic stress; $r \in [0,1]$ is an interpolation parameter~: cases $r\in (0,1)$ are often referred as Jeffreys models whereas case $r=1$ is referred as the Maxwell model. Besides, $\mathbb D(\u)$ (resp. $\W(\u)$) denotes the symmetric  (resp. skew-symmetric) part of the velocity gradient. The function $g_{a}$, $-1\leq a \leq 1$, is a bilinear mapping related to the total derivative, in which the parameter $a$ interpolates between the so-called upper-convected model ($a=1$) and lower-convected model ($a=-1$). Note that the case $a=0$ is known as the so-called corotational model. The function~$g_a$ is defined as
$$
g_{a}(\nabla \u , \bsigma) = \W(\u)\cdot \bsigma - \bsigma \cdot \W(\u) + a(\D(\u)\cdot \bsigma + \bsigma \cdot \D(\u)).
$$

\begin{rem}
In standard derivations of Oldroyd model from kinetic models for dilute polymers, the diffusive term~$\Di \Delta \bsigma$ is routinely omitted, on the grounds that it is several orders of magnitude smaller than the other terms in the equation.
It physically corresponds to a centre-of-mass diffusion term in the dumbell models and it is in the range of about~$10^{-9}$ to~$10^{-7}$ when the macroscopic length-scale of the domain is of order~$1$, see~\cite{BS07}.
\end{rem}

Let us first discuss the transient case. The standard Oldroyd model (without diffusion: $D=0$) has been the subject of intensive studies~\cite{GS90,LM00,MT04}: P.-L.~Lions \& N.~Masmoudi \cite{LM00} proved a global existence result of weak solutions for any data in the corotational case only.  In the general case, C.~Guillop\'e \& J.-C. Saut \cite{GS90} proved the existence and uniqueness of local strong solutions; besides, if the fluid is not too elastic and if the data are sufficiently small, then solutions are global. Then L.~Molinet \& R.~Talhouk \cite{MT04} proved that the smallness assumption on the elasticity of the fluid could be relaxed. The diffusive Oldroyd model (with diffusion: $D>0$) has been studied by a few authors: recently P.~Constantin \& M.~Kliegl \cite{CK12} proved the existence of global strong solutions in 2D for the Cauchy problem and uniqueness of the solution among a class of strong solutions. Notice also that other regularizations of the standard Oldroyd model have been studied, see in particular \cite{BS07}.

Let us now discuss the stationary case. To our knowledge, the only available result is due to M.~Renardy \cite{Re85} and focuses on the standard model (without diffusion) only: existence and uniqueness of strong solutions is proved under the assumption of small regular data. The method used by M.~Renardy is based on a reformulation of the Oldroyd model as a ``Newtonian generalized'' fluid: the contribution of the stress $\div \bsigma$ is expressed as an implicit function of the velocity field~$\u$ and, then, an iterative scheme is built upon this fully nonlinear system. 

Now let us discuss the diffusive Oldroyd model in the stationary case. The mathematical analysis of the model can be approached with a completely different framework, as the diffusive contribution in the stress equation drastically changes the mathematical properties of the system. So far, the method that we present in the present paper is based on a classical weak formulation and then on energy estimates. 

Let us underline the main differences with the standard Oldroyd model studied by M. Renardy:
\begin{itemize}
\item the diffusive model makes it possible to handle with irregular data;
\item in the corotational case, the smallness of the data is not needed anymore.
\end{itemize}

The present paper is composed of five sections. In Section~\ref{sec:results} we give the weak formulation of the problem, the mathematical framework and the main results. The last three sections are devoted to the proof of the main result (Theorem~\ref{theorem}): in Section~\ref{sec:existence} the existence result, in Section~\ref{sec:uniqueness}, the uniqueness result and, in Section~\ref{sec:regularity}, the regularity result.
 
\section{Weak formulation and main results}\label{sec:results}

The variational formulation of Problem~\eqref{problem-strong} is written
\begin{equation}\label{problem-weak}
\left\{
\begin{aligned}
& \textnormal{Find $(\u,\bsigma)\in V\times W$ such that, for all $(\v,\btau)\in V\times W$,}\\
& \Re \int_\Omega ( \u\cdot \nabla \u ) \cdot \v + (1-r) \int_\Omega \nabla \u : \nabla \v + \int_\Omega \bsigma : \D(\v)= \langle \f ,  \v \rangle,\\
& \We \int_\Omega ( \u \cdot \nabla \bsigma + g_a(\nabla \u , \bsigma) ) : \btau + \int_\Omega \bsigma:\btau + \Di \int_\Omega \nabla \bsigma : \nabla \btau \\
& \hspace{8cm} = 2r\dsp \int_\Omega \D(\u) : \btau,
\end{aligned}
\right.
\end{equation}
where the spaces~$V$ and~$W$ are defined by
\begin{equation*}
\begin{aligned}
& V=\{ \u \in H^1_0(\Omega)^3;~\div \u = 0 \},\\
& W=\{\bsigma \in H^1(\Omega)^{3 \times 3};~ \bsigma=\,^T\!\bsigma\},
\end{aligned}
\end{equation*}
and where $\langle \cdot ,  \cdot \rangle$ denotes the duality bracket between~$H^{-1}$ and~$H^1_0$. We will use the following norms:
$$\|\u\|_V^2 = \int_\Omega \|\nabla \u\|^2, \qquad
\|\bsigma\|_W^2 = \int_\Omega \|\bsigma\|^2 + \|\nabla \bsigma\|^2.$$
The main theorem which is proved in this article concerns an existence result for Problem~\eqref{problem-weak}. In its general form, it requires some assumptions on the data. For this, we introduce the following constants:
\begin{equation*}
\begin{aligned}
& C_{(\mathrm{I})} := \frac{8|a| C_\Omega^2 \We \|\f\|_{H^{-1}}}{\min(1-r,\Di)^2},\\[0.1cm]
& C_{(\mathrm{II})} := \frac{\sqrt{2r}\min(1-r,\Di)}{4|a| C_\Omega^2 \We} \left( 1 - \sqrt{1 - C_{(\mathrm{I})}} \right),
\end{aligned}
\end{equation*}
where $C_{\Omega}$ is a constant which only depends on the domain~$\Omega$.

\begin{rem}
\begin{enumerate}
\item[]
\item Note that $C_{(\mathrm{II})}$ is defined provided $C_{(\mathrm{I})}\leq 1$.
\item For $a=0$, we have by continuity $C_{(\mathrm{I})}=0$ and $\dsp C_{(\mathrm{II})}=\frac{\sqrt{2r}\|\f\|_{H^{-1}}}{\min(1-r,\Di)}$.
\item Constant~$C_\Omega$ is related to the Sobolev injection $W\subset L^4(\Omega)^{3\times 3}$: 
$$
\forall \bsigma \in W,\ \|\bsigma\|_{L^4} \leq C_\Omega \|\bsigma\|_W.
$$
\item The only physical parameter that is not involved in the definition of constant~$C_{(\mathrm{I})}$ is the Reynolds number~$\Re$.
In other words, the existence result that is further described does not depend on the value of the Reynolds number.
\end{enumerate}
\end{rem}

\begin{thm}\label{theorem}
Let $\Omega$ be a Lipschitz bounded open set in $\R^3$ and $\f\in H^{-1}(\Omega)^3$.
Let $\Re \geq 0$, $\We\geq 0$, $0< r<1$, $-1\leq a\leq 1$ and $\Di>0$. 
\begin{itemize}
\item {\it Existence.} If $C_{(\mathrm{I})} \leq 1$ then, there exists a solution $(\u, \bsigma)$ of problem~\eqref{problem-weak} which satisfies
\begin{equation}\label{velocity-petitesse3}
2r \|\u\|_V^2 + \|\bsigma\|_W^2 \leq C_{(\mathrm{II})}^2.
\end{equation}
Moreover there exists~$p\in L^2(\Omega)$ such that~$(\u,p)$ satisfies $(\mathrm P_{i})_{i=1,2,3}$ in the sense of distributions.
\item {\it Uniqueness.} Problem~\eqref{problem-weak} admits at most one solution if {\em one} of the following conditions is satisfied:
\begin{enumerate}
\item[a)]~$\|\f\|_{H^{-1}}$ is small enough;
\item[b)]~$\Re$ and~$\We$ are small enough.
\end{enumerate}
\item {\it Regularity.} If $\Omega$ is of class $\mathcal C^\infty$ and if each component of $\f$ belongs to $\mathcal C^\infty(\overline \Omega)$ then each component of any solution of~\eqref{problem-weak} belongs to~$\mathcal C^\infty(\overline \Omega)$ and the considered solution satisfies \eqref{problem-strong} in a classical sense.
\end{itemize}
\end{thm}

\noindent Let us mention two corollaries.

\begin{cor}
Problem~\eqref{problem-weak} admits a unique solution if {\em one} of the following conditions is satisfied:
\begin{enumerate}
\item[a)]~$\|\f\|_{H^{-1}}$ is small enough;
\item[b)]~$\Re$ and~$\We$ are small enough.
\end{enumerate}
\end{cor}

\begin{cor}
If $a=0$, there exists a solution for all data.
\end{cor}

The proof is decomposed into four parts. In section~3 we show the existence for the weak formulation~\eqref{problem-weak} using a Galerkin approximations and compactness results to perform the limit. The existence of a pressure is obtained by De Rham theory. In section 4, we prove the uniqueness of the solution and, in section 5, we investigate the regularity of the weak solutions and prove that, if the data are regular, so is the solution.

\section{Existence result}\label{sec:existence}

As a preliminary, $V\times W$ is endowed with the scalar product $(\cdot,\cdot)_{V\times W}$ defined~by
$$((\u,\bsigma),(\v,\btau))_{V\times W} = 2r(\u,\v)_V + (\bsigma,\btau)_W.$$
As $V$ and $W$ are separable Hilbert spaces, we consider a countable orthonormal basis $(\w_k)_{k\in \N}$ in the space~$V$, and a countable orthonormal basis $(\s_k)_{k\in \N}$ in the space~$W$. We use the notation $V_{k}:=\mathrm{span}(\w_1,...,\w_k)$ and $W_{k}:=\mathrm{span}(\s_1,...,\s_k)$. For each fixed integer $k\in \N$, we would like to define an approximate solution $(\u_k,\bsigma_k)$ of~\eqref{problem-weak} by
$$
\u_k = \sum_{i=0}^k \alpha_{i,k} \w_i, \qquad \bsigma_k = \sum_{i=0}^k \beta_{i,k} \s_i,
$$
satisfying the variational problem
\begin{equation}\label{Galerkin1}
\begin{aligned}
& \Re \int_\Omega ( \u_k\cdot \nabla \u_k ) \cdot \v_k + (1-r) \int_\Omega \nabla \u_k : \nabla \v_k + \int_\Omega \bsigma_k : \D(\v_k)= \langle \f ,  \v_k \rangle,\\
& \We \int_\Omega ( \u_k \cdot \nabla \bsigma_k + g_a(\nabla \u_k , \bsigma_k) ) : \btau_k + \int_\Omega \bsigma_k:\btau_k + \Di \int_\Omega \nabla \bsigma_k : \nabla \btau_k \\
& \hspace{9.1cm} = 2r \int_\Omega \btau_k : \D(\u_{k}),
\end{aligned}
\end{equation}
for all $(\v_{k},\btau_{k})\in V_{k}\times W_{k}$. Equations~\eqref{Galerkin1} form a system of nonlinear equations for $\alpha_{1,k},...,\alpha_{k,k}$, $\beta_{1,k},...,\beta_{k,k}$, and the existence of a solution of this system is not obvious. We use the following lemma:
\begin{lem}\label{lemmTemam}
Let $X$ be a finite dimensional Hilbert space with scalar product $(\cdot,\cdot)$ and norm~$\|\cdot\|$, and let~$P$ be a continuous mapping from $X$ into itself such that
\begin{equation}\label{condi-lemma}
\exists R>0;~ \forall \xi\in X \quad (\|\xi\| = R \Longrightarrow (P(\xi),\xi) \geq 0).
\end{equation}
Then there exists $\xi \in X$, $\|\xi\| \leq R$, such that $P(\xi)=0$.
\end{lem}
The related proof, based on the Brouwer fixed point theorem, can be found in~\cite[p.166]{Temam77}. We only note that the result which is proved in~\cite{Temam77} corresponds to the case where the inegality $(P(\xi),\xi) \geq 0$ in the assertion~\eqref{condi-lemma} is a strict inequality. The case of a large inequality holds too, with the same proof.\\

\noindent
We apply this lemma to prove the existence of $(\u_k,\bsigma_k)$ as follows: let~$X$ be the space defined as $X:=V_{k}\times W_{k}$, endowed with the scalar product inherited from $V\times W$.
Let~$P_k$ the mapping from $X$ into itself defined by, for all $((\u,\bsigma),(\v,\btau))\in X^2$,
\begin{equation*}
\begin{aligned}
& (P_k(\u,\bsigma),(\v,\btau))_X \\
& \hspace{1cm} = 2r \left[ \Re \int_\Omega ( \u\cdot \nabla \u ) \cdot \v + (1-r) \int_\Omega \nabla \u : \nabla \v + \int_\Omega \bsigma : \D(\v) - \langle \f ,  \v \rangle \right] \\
& \hspace{1.5cm} + \left[ \int_\Omega \bsigma:\btau + \Di \int_\Omega \nabla \bsigma : \nabla \btau + \We \int_\Omega ( \u \cdot \nabla \bsigma + g_a(\nabla \u , \bsigma) ) : \btau \right.\\
& \hspace{10.2cm} - 2r \left. \int_\Omega \D(\u) : \btau \right].
\end{aligned}
\end{equation*}
The continuity of~$P_k$ is obvious. Let us show that condition~\eqref{condi-lemma} holds. Denoting $\xi=(\u,\bsigma)\in X$ we have
\begin{equation*}
\begin{aligned}
& (P_k(\xi),\xi)_X = 2r(1-r) \int_\Omega \| \nabla \u \|^2 + \int_\Omega \|\bsigma\|^2 + \Di \int_\Omega \|\nabla \bsigma\|^2 \\
& \hspace{6cm} + \We \int_\Omega g_a(\nabla \u , \bsigma) : \bsigma- 2r \langle \f ,  \u \rangle.
\end{aligned}
\end{equation*}
From the H\"older inequality and the Sobolev injection $W\subset L^4(\Omega)^9$, there exists a positive constant~$C_\Omega$ which only depends on the domain, such that $\|\bsigma\|_{L^4} \leq C_\Omega \|\bsigma\|_W$. Thus, the contribution~$g_a$ satisfied
\begin{equation*}
\int_\Omega g_a(\nabla \u , \bsigma) : \bsigma \leq 2|a| \int_\Omega \| \nabla \u \| \| \bsigma \|^2 \leq 2|a| \| \nabla \u \|_{L^2} \| \bsigma \|_{L^4}^2,
\end{equation*}
which yields
\begin{equation}\label{est1}
\int_\Omega g_a(\nabla \u , \bsigma) : \bsigma \leq 2|a|C_\Omega^2 \|\u \|_V \|\bsigma\|_W^2.
\end{equation}
We deduce that
\begin{equation*}
\begin{aligned}
(P_k(\xi),\xi)_X \geq  & \, 2r(1-r) \|\u \|_V^2 + \min(1,\Di) \|\bsigma\|_W^2 \\[0.1cm]
& \hspace{0cm} - 2|a| C_\Omega^2 \We \|\u \|_V \|\bsigma\|_W^2 - 2r \|\f\|_{H^{-1}} \|\u\|_V.
\end{aligned}
\end{equation*}
Using the definition of the norm $\|\xi\|_X^2 = 2r\|\u\|_V^2 + \|\bsigma\|_W^2$, we obtain
$$\begin{array}{r}
(P_k(\xi),\xi)_X \geq \min(1-r,\Di) \|\xi \|_X^2 - \dsp\frac{\sqrt{2}|a|C_\Omega^2 \We}{\sqrt{r}}\|\xi\|_X^3 - \sqrt{2r} \|\f\|_{H^{-1}} \|\xi\|_X.
\end{array}$$
We write this inequality as $(P_k(\xi),\xi)_X \geq \|\xi \|_X ( -\alpha \|\xi\|_X^2 + \beta \|\xi\|_X - \gamma )$.\\

$\bullet$
If $a\neq 0$ ({\it i.e.} $\alpha>0$), we deduce that $(P_k(\xi),\xi)_X$ may be positive for some value of $\xi$ if the discriminant $\beta^2-4\alpha\gamma$ is nonnegative; in this case we have
$$
(P_k(\xi),\xi)_X \geq 0\quad \Longleftrightarrow \quad 
 \frac{\beta-\sqrt{\beta^2-4\alpha\gamma}}{2\alpha} \leq \|\xi\|_X \leq \frac{\beta+\sqrt{\beta^2-4\alpha\gamma}}{2\alpha}.
$$
In particular, by denoting
$$
C_{(\mathrm{II})}=\frac{\beta-\sqrt{\beta^2-4\alpha\gamma}}{2\alpha},$$
we obtain
$$
\|\xi\|_X=C_{(\mathrm{II})}\quad \Longrightarrow \quad (P_k(\xi),\xi)_X \geq 0.
$$
We also note that $C_{(\mathrm{II})}>0$ (except for the case $\f=\0$ ({\it i.e.} $\gamma=0$) where we can use instead of~$C_{(\mathrm{II})}$ any constant~$C$ such that $0 < C < \frac{\beta}{\alpha}$).\\

$\bullet$
If $a=0$ ({\it i.e.} $\alpha=0$), we easily have
$$
(P_k(\xi),\xi)_X \geq 0\quad \Longleftrightarrow \quad 
 \|\xi\|_X \geq \frac{\gamma}{\beta}.
$$
We then use
$$
C_{(\mathrm{II})}=\frac{\gamma}{\beta}.
$$
As for the case $a\neq 0$ we have $C_{(\mathrm{II})}>0$ (except for $\f=\0$ ({\it i.e.} $\gamma=0$) where we can use instead of~$C_{(\mathrm{II})}$ any positive constant).\\

\noindent
Finally, assumptions given by Lemma~\ref{lemmTemam} are satisfied taking $R=C_{(\mathrm{II})}$ as soon as $\beta^2-4\alpha\gamma\geq 0$, which is equivalent to the condition~$C_{(\mathrm{I})} \leq 1$ introduced in Theorem~\ref{theorem}.
By Lemma~\ref{lemmTemam}, we deduce that, for any $k\in \N$, Equations~\eqref{Galerkin1} admit a solution~$(\u_k,\bsigma_k)$ which satisfies
\begin{equation}\label{velocityk-petitesse}
2r \|\u_k\|_V^2 + \|\bsigma_k\|_W^2 \leq C_{(\mathrm{II})}^2.
\end{equation}
This estimate~\eqref{velocityk-petitesse} implies that the sequence $(\u_k,\bsigma_k)$ remains bounded in $V\times W$. Thus there exists some $(\u,\bsigma)\in V\times W$ and a subsequence (still denoted by~$k$) such that $(\u_k, \bsigma_k) \rightharpoonup (\u,\bsigma)$ for the weak topology of $V\times W$, as $k \rightarrow +\infty$. As $H^1 \subset L^2$ with compact injection, $(\u_k, \bsigma_k) \to (\u,\bsigma)$ for the strong topology of $L^2(\Omega)^3\times L^2(\Omega)^9$. The consequences are twofold:
\par\noindent
\begin{itemize}
\item we can pass to the limit in all the terms of Equations~\eqref{Galerkin1}, and deduce that $(\u,\bsigma)$ is a solution of~\eqref{problem-weak}~;
\item due to the usual property of the weak limit, we have
$$2r \|\u\|_V^2 + \|\bsigma\|_W^2 \leq\displaystyle \liminf_{k\rightarrow +\infty} 2r \|\u_k\|_V^2 + \|\bsigma_k\|_W^2\leq C_{(\mathrm{II})}^2.$$
\end{itemize}
This concludes the proof of existence of a solution $(\u, \bsigma)$ of problem~\eqref{problem-weak} satisfying Eq.~\eqref{velocity-petitesse3}.\\

\noindent
Now let us prove the existence of a pressure field associated to the incompressibility condition. For a solution~$(\u,\bsigma)$ of problem~\eqref{problem-weak}, we have, for all $\v\in \mathcal D(\Omega)$ such that $\div \v = 0$,
$$\langle \Re (\u\cdot \nabla \u) - (1-r)\Delta \u - \div \bsigma - \f, \v \rangle = 0.$$
By De Rham theorem, there exists a pressure $p\in \mathcal D'(\Omega)$ such that Eq.~\eqref{problem-strong} holds in~$\mathcal D'(\Omega)$. The regularity of $\Delta \u$, $\u\cdot \nabla \u$, $\div \bsigma$ and $\f$ implies that the pressure~$p$ is more regular. For instance, by Sobolev embeddings, we have $$\u\in V \subset L^6(\Omega)^3$$
and since $\nabla u \in L^2(\Omega)^{3\times 3}$, then the convective term $\u\cdot \nabla \u$ belongs to~$L^{3/2}(\Omega)$. Thus, we have
$$
\nabla p=\Delta \u+\div \bsigma +\f -\mathrm{Re} (\u \cdot \nabla u) \in H^{-1}(\Omega).
$$
This implies that~$p \in L^{2}(\Omega)$, see \cite[p.14]{Temam77}.

\section{Uniqueness of the solution}\label{sec:uniqueness}

Let $(\u_1,\bsigma_1)$ and~$(\u_2,\bsigma_2)$ be two solutions of~\eqref{problem-weak} and introduce the difference $(\u,\bsigma) = (\u_2,\bsigma_2) - (\u_1,\bsigma_1)$.
By subtraction we obtain, for all $(\v,\btau)\in V\times W$,
\begin{equation*}
\begin{aligned}
& \Re \int_\Omega ( \u\cdot \nabla \u_2 + \u_1\cdot \nabla \u ) \cdot \v + (1-r) \int_\Omega \nabla \u : \nabla \v + \int_\Omega \bsigma : \D(\v) = 0,\\
& \We \int_\Omega ( \u \cdot \nabla \bsigma_2 + \u_1 \cdot \nabla \bsigma + g_a(\nabla \u , \bsigma_2) + g_a(\nabla \u_1 , \bsigma)) : \btau + \int_\Omega \bsigma:\btau\\ 
& \hspace{6.5cm} + \Di \int_\Omega \nabla \bsigma : \nabla \btau =\displaystyle 2r \int_\Omega \btau : \D(\u).
\end{aligned}
\end{equation*}
Taking $(\v,\btau) = (2r\, \u,\bsigma)$ and adding the two equations, we obtain
\begin{equation*}
\begin{aligned}
& 2r (1-r) \|\u\|_V^2 + \min(1,\Di) \|\bsigma\|_W^2 \leq -2r \Re \int_\Omega ( \u\cdot \nabla \u_2 ) \cdot \u \\
& \hspace{4cm} - \We \int_\Omega ( \u \cdot \nabla \bsigma_2 + g_a(\nabla \u , \bsigma_2) + g_a(\nabla \u_1 , \bsigma)) : \bsigma.
\end{aligned}
\end{equation*}
By definition of the bilinear function~$g_a$, we have
$$
\int_\Omega g_a(\nabla \u , \bsigma_2): \bsigma = (a+1) \int_\Omega (\nabla \u \cdot \bsigma_2): \bsigma + (a-1)\int_\Omega (\bsigma_2 \cdot \nabla \u): \bsigma.
$$
Using the H\"older inequality and the Sobolev injection $W\subset L^4(\Omega)^9$ we deduce that
$$
\left| \int_\Omega g_a(\nabla \u , \bsigma_2): \bsigma \right| \leq 2 C_\Omega^2 \|\u\|_V \|\bsigma_2\|_W \|\bsigma\|_W.
$$
Estimating the other terms as in Equation~\eqref{est1}, we obtain
\begin{equation*}
\begin{aligned}
2r (1-r) \|\u\|_V^2 + \min(1,\Di) \|\bsigma\|_W^2 \leq & \, 2rC_\Omega^2 \Re \| \u_2 \|_V \|\u\|_V^2 \\[0.1cm]
& + 3 C_\Omega^2 \We \|\bsigma_2\|_W \|\u\|_V \|\bsigma\|_W \\[0.1cm]
& + 2|a| C_\Omega^2 \We \|\u_1\|_V \|\bsigma\|_W^2.
\end{aligned}
\end{equation*}
As the solutions~$(\u_1,\bsigma_1)$ and~$(\u_2,\bsigma_2)$ satisfy~Equation~\eqref{velocity-petitesse3}, we deduce that
\begin{equation}\label{estimate18}
\begin{aligned}
2r (1-r) \|\u\|_V^2 + \min(1,\Di) \|\bsigma\|_W^2 \leq \frac{C_\Omega^2 C_{(\mathrm{II})}}{\sqrt{2r}} \Big( & 2 r \Re \|\u\|_V^2\\
& + 3 \We \sqrt{2r} \|\u\|_V \|\bsigma\|_W \\
& + 2 |a| \We \|\bsigma\|_W^2 \Big).
\end{aligned}
\end{equation}
From Young inequality we deduce
$$
3 \We \sqrt{2r} \|\u\|_V \|\bsigma\|_W \leq 2r \Re \|\u\|_V^2 + \frac{9\We^2}{4\Re}\|\bsigma\|_W^2.
$$
Consequently, Estimate~\eqref{estimate18} reads $2rA\|\u\|_V^2 + B \|\bsigma\|_W^2 \leq 0$ where $A$ and $B$ are given by
\begin{equation*}
\begin{aligned}
& A= (1-r) - \frac{4 \Re C_\Omega^2 C_{(\mathrm{II})}}{\sqrt{2r}},\\
& B= \min(1,\Di) - \frac{C_\Omega^2 C_{(\mathrm{II})}}{\sqrt{2r}} \left( \frac{9\We^2}{4\Re} + 2|a|\We \right).
\end{aligned}
\end{equation*}
Using the fact that $C_{(\mathrm{II})}$ tends to~$0$ as $\|\f\|_{H^{-1}}$ tends to~$0$, we deduce that for $\|\f\|_{H^{-1}}$ small enough, the coefficients~$A$ and~$B$ are positive.
In the same way, if~$\Re$ is small enough then $A>0$ and, if~$\We$ is small, then $B>0$.
These inequalities imply $\|\u\|_V=\|\bsigma\|_W=0$, which means $(\u_1,\bsigma_1)=(\u_2,\bsigma_2)$.

\section{Strong solutions}\label{sec:regularity}

As previously, we notice that $\u\in  L^6(\Omega)^3$, so that  $\u\cdot \nabla \u \in L^{3/2}(\Omega)^3$. We also have $\div \bsigma\in L^2(\Omega)^3$ and if $\f$ is regular then the regularity of the Stokes problem, see \cite{Ca61}, implies that 
$$\u\in W^{2,3/2}(\Omega)^3,\qquad p\in W^{1,3/2}(\Omega).$$ 
In the same way, $\u\cdot \nabla \bsigma$ and $g_a(\nabla \u,\bsigma)$ belong to~$L^{3/2}(\Omega)^{3\times 3}$. The regularity of the Laplace problem implies that 
$$\bsigma \in W^{2,3/2}(\Omega)^{3\times 3}.$$ 
Repeating such a process, we find that $\u\in W^{s,3/2}(\Omega)^3$, $p\in W^{s-1,3/2}(\Omega)$ and~$\bsigma \in W^{s,3/2}(\Omega)^9$ for any $s\in \N$. The proof is concluded.


\end{document}